\documentclass[11pt]{article}

\usepackage{amsmath, amsfonts, amssymb, amsthm}

\usepackage{natbib}

\usepackage{setspace}

\DeclareMathSymbol{\leqslant}{\mathalpha}{AMSa}{"36} 

\DeclareMathSymbol{\geqslant}{\mathalpha}{AMSa}{"3E} 

\DeclareMathSymbol{\eset}{\mathalpha}{AMSb}{"3F}     

\newcommand{\ban}{\begin{align}}
\newcommand{\ean}{\end{align}}

\newcommand{\ba}{\begin{align*}}
\newcommand{\ea}{\end{align*}}

\newcommand{\be}{\begin{eqnarray*}}
\newcommand{\ee}{\end{eqnarray*}}

\newcommand{\ben}{\begin{eqnarray}}
\newcommand{\een}{\end{eqnarray}}

\theoremstyle{plain}
\newtheorem{theo}{Theorem}[section]

\newtheorem{lemma}[theo]{Lemma}

\theoremstyle{definition}

\newtheorem{remark}[theo]{Remark}

\begin{document}

\vglue20pt \centerline{\huge\bf Dispersion of volume under the action of }

\medskip

\centerline{\huge\bf  }

\medskip

\centerline{\huge\bf  isotropic  Brownian flows}
 \medskip


\bigskip

\bigskip

\centerline{by}

\bigskip

\medskip

\centerline{{\Large  G.\ Dimitroff\footnotemark[1] and  M.\ Scheutzow\footnotemark[2]}}

\footnotetext[1]{Fraunhofer ITWM, Fraunhofer-Platz 1, D-67663 Kaiserslautern }

\footnotetext[2]{Institut f\"ur Mathematik, MA 7-5, Technische Universit\"at Berlin, Stra\ss e des 17. Juni 136, D-10623 Berlin }

\bigskip

\bigskip

{\leftskip=1truecm

\rightskip=1truecm

\baselineskip=15pt

\small

\noindent{\slshape\bfseries Summary.} We study transport properties of isotropic Brownian flows. Under a transience condition for
the two-point motion, we show asymptotic normality of the image of a finite measure under the flow and -- under slightly stronger
assumptions -- asymptotic normality of the distribution of the volume of the image of a set under the flow. Finally, we show that
for a class of isotropic flows, the volume of the image of a nonempty open set (which is a martingale) converges to a random variable which is
almost surely strictly positive. \bigskip


\bigskip


}


\newcommand{\oq}{{\langle}}
\newcommand{\ve}{{\varepsilon}}
\newcommand{\cq}{{\rangle}_t}
\newcommand{\dd}{{\mathrm{d}}}
\newcommand{\Dd}{{\mathrm{D}}}
\newcommand{\cd}{{\cdot}}
\newcommand{\nix}{{\varnothing}}
\newcommand{\N}{{\mathbb N}}
\newcommand{\Z}{{\mathbb Z}}
\newcommand{\R}{{\mathbb R}}
\newcommand{\Q}{{\mathbb Q}}
\newcommand{\E}{{\mathbb E}}
\renewcommand{\P}{{\mathbb P}}
\newcommand{\F}{{\cal F}}
\newcommand{\C}{{\mathbb C}}
\newcommand{\K}{{\mathbb K}}
\newcommand{\B}{{\cal B}}
\newcommand{\G}{{\cal G}}
\newcommand{\D}{{\cal D}}
\newcommand{\Zz}{{\cal Z}}
\newcommand{\Ll}{{\cal L}}
\newcommand{\A}{{\cal A}}
\newcommand{\Po}{{\cal P}}
\newcommand{\Sy}{{\cal S}}
\newcommand{\cZ}{{\cal Z}}
\newcommand{\M}{{\cal M}}
\newcommand{\Nn}{{\cal N}}
\newcommand{\p}{{\mathbf P}}
\newcommand{\X}{{\mathbb X}}
\newcommand{\interior}[1]{{\mathaccent 23 #1}}
\newcommand{\bem}{\begin{em}}
\newcommand{\eem}{\end{em}}
\def\eins{{\mathchoice {1\mskip-4mu\mathrm l}
{1\mskip-4mu\mathrm l}
{1\mskip-4.5mu\mathrm l} {1\mskip-5mu\mathrm l}}}
\newcommand{\olim}[1]{\begin{array}{c} ~\\[-1.4ex] \overline{\lim} \\[-1.35ex]
                       {\scriptstyle #1}\end{array}}
\newcommand{\plim}[2]{\begin{array}{c} #1\\[-1.175ex] \longrightarrow \\[-1.2ex]
                       {\scriptstyle #2}\end{array}}

\newcommand{\ulim}[1]{\begin{array}{c} ~\\[-1.175ex] \underline{\lim} \\[-1.2ex]
                       {\scriptstyle #1}\end{array}}

\renewcommand{\theequation}{\thesection.\arabic{equation}}

\section{Introduction}\label{intro}
It has been suggested that stochastic flows can be used as a model for studying the spread of passive tracers within a
turbulent fluid. Individual particles perform diffusions while the motions of adjacent particles are correlated
and form a stochastic flow of homeomorphisms $\{\phi_{s,t}\}_{0 \le s \le t}$.
We consider isotropic Brownian flows (IBF) on $\R^d$ -- a special class of stochastic flows characterized by their spatial translation and rotational invariance, temporal homogeneity and independence of their increments. They have been
extensively studied by many authors e.g.  Baxendale and Harris (\cite{BaH86}, \cite{Bax91}), Le Jan (\cite{LeJ85}, \cite{LeJ87},
\cite{DaLe88}), Yaglom \cite{Yag57}, Cranston, Scheutzow and Steinsaltz \cite{CSS99}. We will recall basic facts about IBFs in
the next section.   \\
The one-point motion of an isotropic Brownian flow is a  Brownian motion, but the distribution of the two-point motion
$(\phi_t(x), \phi_t(y))^T$ is  not Gaussian (we  write $\phi_t$ for $\phi_{0,t}$).   In the transient case, in which the distance
$|\phi_t(x)-\phi_t(y)|$ grows to infinity as  $ t\to \infty $, the correlation between $\phi_t(x)$ and $\phi_t(y)$ decays to zero
and it is reasonable to expect some kind of joint asymptotic normality of the empirical distribution of $\left\{\phi_t (x), x \in B\right\}$ for
suitable initial sets $B$. In the following, we fix $d \ge 2$ and let $\lambda$ be Lebesgue measure on $\R^d$.
We will consider the following two questions:
\begin{enumerate}
\item[1.] Given two ``nice" sets $A,\,B\subset \R^d$, what portion of the set $B$ moves under the action of the flow at time $t$
into the set $\sqrt t A$? That is, how does
\begin{align*}
\frac{1}{\lambda(B)}\lambda\left (  B\cap \phi_t^{-1}(\sqrt t A)\right)
\end{align*}
behave as $t \to \infty$?
\item[2.]Given two ``nice" sets $A,\,B\subset \R^d$, what portion of the set $\phi_t(B)$ is in the set $\sqrt t A$? That is, how does
\begin{align*}
\frac{1}{\lambda(\phi_t(B))}\lambda\left ( \phi_t(B)\cap \sqrt t A \right)
\end{align*}
behave as $t \to \infty$?
\end{enumerate}
Theorems \ref{an_main1} and \ref{an_main2} below  deal with these questions.
The asymptotic normality result arising from the second question is however weaker in the sense that it holds only on the set
where  $\lim_{t\to \infty}\lambda(\phi_t(B))> 0$. In case the top Lyapunov exponent of the flow is strictly positive,
results of Le Jan on the statistical equilibrium  show that
$\big(\lambda(\phi_t(B))\big)_t$ is a uniformly integrable martingale and thus its limit cannot be trivial, i.e.
$\P\big[\lim_{t\to \infty}\lambda(\phi_t(B))> 0\big]>0$ whenever $\lambda(B) > 0$. One can ask, if the limit is almost surely
strictly positive, that  is if the volume persists almost surely. In Section 4 we answer this question affirmatively
for a class of IBFs characterized by a simple condition on the parameters. \\
In \cite{DKK04}, Question 1.~is  considered  for rather  general measures and   for stochastic flows driven by a finite
dimensional Brownian motion on a compact manifold (or periodic flows on $ \R^d $). Using different techniques
relying crucially on the compactness of the state space, they show asymptotic normality.\\
In \cite{CiZi97} and \cite{Zi97}, the authors treat  the asymptotic behavior of the first and the second moments of
$ m\circ \phi_t^{-1} $ for some finite measure $ m $ on $ \R^d $. This is a somewhat different approach for describing
the dispersion of sets. Observe that  Question 1 above can be formulated as the asymptotic behavior of
$ \lambda_B \circ \psi_t^{-1} $ with $ \psi_t:=\frac{\phi_t}{\sqrt t} $, where
$ \lambda_B(\dd x)=\frac{1}{\lambda(B)}\eins_B(x)\lambda(\dd x) $. This suggests to generalize Question 1 by
considering the transport of an arbitrary finite measure $m$ rather than just $\lambda_B$.\\

\section{Isotropic Brownian flows}
We provide a brief introduction to IBFs following mainly \cite{BaH86}, \cite{LeJ85} and \cite{Yag57}. A rather extensive
collection of properties can be found in \cite{Dim06}. \\

A \bem (forward) stochastic flow of homeomorphisms \eem on $\R^d$ is a family of random homeomorphisms  $\left\{\phi_{s,t} : 0\le s\le  t< \infty\right\}$ of $\R^d$ into itself, such that (up to a set of measure zero)  $\phi_{u,t}\circ \phi_{s,u}=\phi_{s,t}$ for $s\le u\le t$ and
$(s,t,x) \mapsto \phi_{s,t}(x)$ is continuous.
The flow is called \bem Brownian \eem if the increments are independent, i.e.~for arbitrary $s_1\le t_1\le s_2\le \dots \le t_n$,
the increment mappings  $\phi_{s_1,t_1},\dots,\phi_{s_n,t_n}$ are independent. If  the increments are stationary, the flow is
called \bem temporally homogeneous. \eem  A flow is called \bem spatially homogeneous \eem if its distribution  is invariant with
respect to translations $T:\R^d\to \R^d$
in the sense that  the flows $T\circ \phi $ and $\phi\circ T$
have the same law. A spatially and temporally homogeneous Brownian flow is called {\em isotropic Brownian flow} (IBF) if,
in addition, $\phi$ is invariant under orthogonal transformations on $\R^d$.\\
For an IBF $\phi$, the one-point motion $t \mapsto \phi_{0,t}$ is a $d-$dimensional Brownian motion. The covariance tensor
$b$ of $\phi$ is defined as
$$
b_{i,j}:=\frac{\dd \langle \phi_{0,t}^i(x),\phi_{0,t}^j(0)\rangle}{\dd t}|_{t=0},
\qquad x \in \R^d,\;i,j \in {1,...,d},
$$
where $ \langle \phi_{0,t}^i(x),\phi_{0,t}^j(0)\rangle$ denotes the covariation process of $\phi_{0,t}^i(x)$ and $\phi_{0,t}^i(0)$.
The properties of $\phi$ then imply that $\dd \langle \phi_{s,t}^i(x),\phi_{s,t}^j(y)\rangle = b_{i,j}(\phi_{s,t}(x)-
\phi_{s,t}(y))\dd t$ for all $x,y \in \R^d,\,0 \le s \le t$. Following \cite{BaH86}, we assume that $b(0)=\mathrm{Id}_{\R^d}$
(which is just a normalization condition) and that $x \mapsto b(x)$ is $C^4$. We will further assume that $d \ge 2$ and that
$\lim_{|x| \to \infty} b(x)=0$ (so the random translation part of $\phi$ is zero).

According to  \cite{BaH86} and \cite{Yag57}, the positive semidefinite function $b$ has necessarily the form
\begin{align}\label{pre_ibf2}
b_{i,j}(x)= \left\{ \begin{array}{ll} \left(B_L(|x|)-B_N(|x|)\right)\frac{x_ix_j}{|x|^2}+\delta_{i,j}B_N(|x|) & \text{for }x\ne 0  \\ [0.4cm]
\delta_{i,j} & \text{for }x= 0\,,
\end{array} \right.
\end{align}
where $B_L$ and $B_N$ are the so called longitudinal and  transversal covariance  functions defined by
\begin{align*}
&B_L(r)=b_{p,p}(re_p), \,\,\, r\ge 0 \,\,\,\text{ and } \\
&B_N(r)=b_{p,p}(re_q), \,\,\, r\ge 0,\,\,\, p\ne q \,,
\end{align*}
where $e_i$ denotes the $ i $-th standard basis vector in $ \R^d $.
$B_L$ and $B_N$ are bounded $C^4$ functions with bounded derivatives. Set
\begin{align*}
\beta_L&:=-\partial_p\partial_pb_{p,p}(0)=-B_L''(0)\\
\beta_N&:=-\partial_q\partial_q b_{p,p}(0)=-B_N''(0) \,\,\, q\ne p\,,
\end{align*}
where $B_L''(0)$ and $B_N''(0)$ denote the right-hand second derivatives of $B_L$ and $B_N$ at zero. The asymptotics of
$B_L$ and $B_N$ around zero are given by
\begin{align}\label{pre_ibf3}
\begin{split}
&B_L(r)=1-\frac{1}{2}\beta_L r^2+O(r^4)\,,\hspace{0.5cm} \text{ for } r\to + 0\,\,\,\text{ and }\\
&B_N(r)=1-\frac{1}{2}\beta_N r^2+O(r^4)\,,\hspace{0.5cm} \text{ for } r\to +0,
\end{split}
\end{align}
where $ \beta_L,\,\beta_N>0 $. Moreover, $ |B_L(r)|,\,|B_N(r)|<1 $ for arbitrary $ r>0 $.
 The constants $\beta_L$ and $\beta_N$ satisfy
\begin{align}\label{betta}
0< \frac{d-1}{d+1}\le \frac{\beta_L}{\beta_N}\le 3\,.
\end{align}

We will abbreviate $\phi_t:=\phi_{0,t}$. According to Theorem 4.2.4 from \cite{Ku90},
the $n$-point motion $\left\{ (\phi_{t}(x_1), \dots,\phi_{t}(x_n))^T : t\ge 0\right\}$
is a diffusion on $\R^{nd}$ with generator
\begin{align*}
\mathcal L_n f(x_1,\dots,x_n)=\sum_{p,g=1}^n\sum_{i,j=1}^d b_{i,j}(x_p-x_q)\frac{\partial^2f}{\partial x_p^i \partial x_q^j}(x_1,\dots,x_n)
\end{align*}
for $f\in C^2_b(\R^{nd}:\R)$.
An important and useful consequence of the  isotropy and the translation invariance is that the distance process
$\rho_t:=\rho_t(x,y):=|\phi_{t}(x)-\phi_t(y)|$ is a diffusion. Its generator is (see Lemma 3.10 and Corollary 3.16 in \cite{BaH86})
\begin{align*}
\mathcal L g(r)=(1-B_L(r))\frac{\dd^2 g}{\dd r^2}(r)+(d-1)\frac{1-B_N(r)}{r}\frac{\dd g}{\dd r}(r) \,\,\,\text{ for all } g \in C^2_b.
\end{align*}
The independence of the increments and the spatial homogeneity imply that  the law of the matrix-valued linearization
$L_t:=D\phi_{t}(x)$ of the flow around the trajectory started in $x\in \R^d$ does not depend on $x$ and
has itself stationary and independent increments. Therefore it follows from Oseledec's multiplicative ergodic theorem (\cite{Arn98}) that
$$\lim_{t\to \infty}\frac{1}{2t}\log s_i\left[L_t^T L_t\right]=:\lambda_i $$
exists almost surely, where $s_i[M]$ denotes the $i$-th largest eigenvalue of the symmetric matrix $M$. Moreover,
$\lambda_1\ge \dots\ge \lambda_d$ are constant and are called \bem Lyapunov exponents\eem. The following formula for
the Lyapunov exponents can be found for example in \cite{BaH86} and \cite{LeJ85}:
$$\lambda_i=(d-i)\frac{\beta_N}{2}-i\frac{\beta_L}{2} \,\,\,\,\,\,\,\,\text{ for }\,\,\,\,\,\, i=1.\dots,d.$$
The recurrence/transience  modes of the distance process  can be characterized by the top Lyapunov exponent
$\lambda_1$ (see \cite{BaH86} and \cite{LeJ85}). In particular we have the following
\begin{remark}\label{ibf_rectrans} For arbitrary $x\ne y$
\begin{align}\label{trans_ibf}
|\phi_t(x)-\phi_t(y)| \underset{t \to +\infty} \longrightarrow
\infty \textrm{    in probability}\,,
\end{align}
if one of the following holds\\
(i) $d\ge 4$ \hspace{1cm}or\hspace{1cm} (ii) $d=3$ and $\lambda_1\ge 0$ \hspace{1cm}or\hspace{1cm} (iii) $d=2$ and $\lambda_1> 0$.\\
We will call an IBF {\em transient} if it satisfies \eqref{trans_ibf}. These conditions are also necessary with the
possible exception of the case $d=2,\,\lambda_1=0$ which seems to be unclear.
  The above convergence  holds also almost surely in the cases (i) and (ii), but not in case (iii).
\end{remark}
As pointed out in \cite{LeJ85}, the  measure $\psi(|x-y|)\lambda(\dd x)\lambda(\dd y)$, where
 \begin{align}\label{ibf2p_inv_meas}
\psi(s)= \frac{1}{1-B_L(s)}\exp \left [ -(d-1)\int\limits_s^\infty\frac{B_L(u)-B_N(u)}{u(1-B_L(u))} \dd u \right ]
\end{align}
is a reversible  invariant measure for the two-point motion (the integral in \eqref{ibf2p_inv_meas} is always finite).
Further we have
\begin{align}\label{ibf2p_inv_meas_asymp}
\lim_{s \to +\infty}\psi(s)=1\,\,\,\text{ and }\,\,\,
\psi(s)\sim \frac{c}{\beta_L}s^{(d-1)\frac{\beta_N}{\beta_L}-(d+1)} \,\, \text{ for } s \to +0 ,
\end{align}
where $ c $ is a strictly positive constant. We conclude this section by mentioning, that an IBF can be
represented as the solution flow of a {\em Kunita-type} stochastic differential equation, see \cite{CSS99}.

\section{The dispersion of volume}
In this section, we will deal with the two questions about asymptotic normality formulated in the introduction.
We will assume throughout this section, that
$\phi$ is a transient IBF. As pointed out at the end of the introduction, we will generalize the first question by
considering the transport of finite measures rather than of sets of finite Lebesgue measure.
\begin{theo}\label{an_main1}
Let $\left \{ \phi_{s,t}(x,\omega) : s,t \in \R_+ \right \} $ be a transient IBF and let $m$ be an atomless nonzero
finite  measure on $(\R^d,\B^d)$.
Further, let the random measure $l_t$ on the Borel sets of $\R^d$ be defined as
$$
l_t(A):=\frac{1}{m(\R^d)}m(\{x \in \R^d:\phi_t(x) \in \sqrt{t}A\}).
$$
Then we have
\begin{displaymath}
l_t(A) \underset{t \to +\infty}\longrightarrow
\mathcal{N}(0, \textnormal{Id}_{\R^d}) (A)\textrm{ in } \mathrm{L}^2\,,\nonumber
\end{displaymath}
for every $A \in  \B(\R^d)$. In particular, $l_t$ converges to $\mathcal{N}(0, \textnormal{Id}_{\R^d})$ weakly in
probability.
\end{theo}

\begin{theo}\label{an_main2}
Let $\left \{ \phi_{s,t}(x,\omega): s,t \in \R_+ \right \}$   be an IBF with $\lambda_1>0$
and  $B\in \B(\R^d)  $ be a Borel set with nontrivial Lebesgue measure  $ \lambda(B)\in (0,\infty) $.
Then, for any set $A \in \B(\R^d)  $, we have
\begin{align} \label{qares1}
\lambda\left( \phi_t(B)\cap \sqrt{t}A \right)-\mathcal{N}(0, \textnormal{Id}_{\R^d})\left(A\right )\lambda\left(\phi_t(B)\right )
\underset{t \to +\infty} \longrightarrow 0 \textrm{ in } \mathrm{L}^2 \textrm{. }
\end{align}
Moreover, conditioned on the set $\Omega_B:=\left \{ \lim\limits_{t \to \infty}\lambda(\phi_t(B))\ne 0\right\}$,
\begin{align} \label{qares2}
\frac{1}{\lambda(\phi_t(B))}\lambda\left(\phi_t(B)\cap \sqrt{t}A \right )\underset{t \to +\infty} \longrightarrow
\mathcal{N}(0,\textnormal{Id}_{\R^d})(A) \text{ in }\P(\cdot\,\,| \,\,\Omega_B)\textrm{  probability}
\end{align}
\end{theo}
\begin{remark}
We will see in Section \ref{persistence}, that the limit
$\lim_{t \to \infty} \lambda(\phi_t(B))$ exists almost surely and is finite and strictly positive with positive probability,
i.e.~$\Omega_B$ has positive probability. It is still an open question, whether under the condition of
 Theorem  \ref{an_main2} we always have $\P\left(\Omega_B\right)=1$, but we will see that under the additional
condition $\frac{\beta_N}{\beta_L}>\frac{d}{d-1}$ we indeed have $\P(\Omega_B)=1$.
\end{remark}

\subsection{Proofs of the dispersion results}
The proofs of  Theorems \ref{an_main1} and \ref{an_main2} rely on the following two simple lemmas:
\begin{lemma}\label{an_clt4mg}
Let $ (\Omega,\F,(\F_t)_{t \in \R_+},\P) $ be a filtered probability
space, satisfying the usual conditions. Further, let $\left\{ M_t : t \in \R_+\right\}$ be a continuous local $\F_t$-martingale
starting at $m\in \R$,  with quadratic variation satisfying
\begin{equation}\label{an_almostbm}
\lim\limits_{t \to +\infty}\frac{1}{t}\langle M\rangle_t=c \textrm{
  in probability }
\end{equation}
for some strictly positive constant c. Then
\begin{equation}\label{an_clt4mg_1}
\frac{1}{\sqrt{t}} M_t \underset{t \to +\infty} \longrightarrow
\mathcal{N}(0,c)\textrm{ in law, }
\end{equation}
where $\mathcal{N}(0,c) $ denotes the zero mean Gaussian distribution
with variance $c$.
\end{lemma}
\noindent \textbf{Proof}: \\
It suffices to consider the case $m=0$ and $c=1$. The continuous local martingale $M$ can be represented as a time-changed
Brownian motion (on a possibly enriched probability space), see  \cite{KaShr88}, p.174-175, i.e.~there exists a standard Brownian
motion $B$ and an increasing family $\tau(t)$, $t \ge 0$ of random times such that $M_t=B(\tau(t))$ for all $t \ge 0$. For $T>0$, define a new Brownian
motion $\tilde B_T(t):=\frac{B(tT)}{\sqrt{T}}$, $t \ge 0$. Due to  (\ref{an_almostbm}), we have $\tau(t)/t \to 1$ in probability and
therefore
$$
\frac{M_T}{\sqrt{T}} = \tilde B_T \left( \frac{\tau(T)}{T} \right) \to \mathcal{N}(0,1) \textrm{ in law, }
$$
and the proof is complete. \hfill $ \square $

\begin{lemma}\label{an_CLTforIBF}
Let $\left \{\phi_{s,t}(x,\omega): s,t \in \R_+ \right \}$  be a transient IBF. For every $x,y \in \R^d$, $x \ne y$ the following convergence holds:
\begin{displaymath}
\frac{1}{\sqrt{t}}\left( \begin{array}{c}\phi_t(x) \\ \phi_t(y)
  \end{array}\right) \underset{t \to +\infty} \longrightarrow
\mathcal{N}(0,\textnormal{Id}_{\R^{2d}}) \textrm{ in distribution. }\nonumber
\end{displaymath}
\end{lemma}

\noindent \textbf{Proof}: \\
The proof relies on the Cram\'er-Wold Theorem (see e.g. \cite{Kall02}). \\
For arbitrary $\xi=(\alpha_1, \ldots,\alpha_d,\beta_1,
\ldots,\beta_d) \in \R^{2d}$ satisfying $|\xi|=1$,
\begin{equation}
\psi_t:=\sum_{i=1}^{d}\alpha_i\phi_t^i(x)+\sum_{i=1}^{d}\beta_i\phi_t^i(y)\nonumber
\end{equation}
is a continuous $\F_t$-martingale with quadratic variation
\begin{align}\label{an_quadr}
\langle \psi \rangle_t = t\sum_{i=1}^{d}\left( \alpha_i^2 + \beta_i^2 \right )+\sum_{i,j=1}^{d}\alpha_i
\beta_j\int\limits_0^t
b_{i,j}(\phi_s(x)-\phi_s(y))\dd s .
\end{align}
Since $\lim_{|x|\to +\infty}b_{i,j}(x)=0$ for all $i,j \in \{1,\ldots,d\}$ and $\phi$ is transient, we obtain
\begin{align*}
\frac{1}{t}\int\limits_0^t b_{i,j}(\phi_s(x)-\phi_s(y))\dd s \underset{t
  \to +\infty}\longrightarrow  0 \hspace{0.5cm }\textrm{in } \mathrm{L}^1,
\end{align*}
and hence
\begin{align*}
\frac{1}{t}\langle \psi \rangle_t \underset{t \to +\infty}\longrightarrow  \sum_{i=1}^{d}\left(
  (\alpha_i)^2 +(\beta_i)^2 \right )=|\xi|^2=1 \textrm {   in }  \mathrm{L}^1.
\end{align*}
$\psi_t$ satisfies  the conditions of Lemma \ref{an_clt4mg} and the statement of the lemma
therefore follows from the Cram\'er-Wold Theorem. \hfill$\square$\\[0.5cm]

Now we are ready to prove the main theorems:\\

\noindent\textbf{Proof of  Theorem \ref{an_main1}}:\\
Without loss of generality, we assume that $m$ is a probability measure on $\R^d$.
Fix a Borel set $A$. Fubini's Theorem and Lebesgue's dominated convergence theorem imply that
\begin{align*}
&\E \,l_t(A)=\E \int\limits_{\R^d}  \eins_{\sqrt{t}A}(\phi_t(x))m(\dd x)=\int\limits_{\R^d} \P(\phi_t(x)\in \sqrt t A)m(\dd x)\\[0.4cm]
& =\int\limits_{\R^d} \P(\phi_1(0) \in A-\frac{x}{\sqrt t})m(\dd x)= \int\limits_{\R^d}
  \,\,\mathcal{N}(0,\text{Id}_{\R^d}) (A-\frac{x}{\sqrt t}) m(\dd x)
\underset{t \to +\infty}\longrightarrow \mathcal{N}(0, \text{Id}_{\R^d}) (A)\textrm{ .}
\end{align*}
Assuming in addition that the set $A$ is closed and using  Lemma \ref{an_CLTforIBF}  and the
fact that  $m\otimes m$ does not charge the diagonal $D:=\{(x,x) : x\in \R^d\}$, we get
\begin{align*}
&\limsup_{t \to \infty}\E( l_t(A))^2 = \limsup_{t \to \infty} \E \int\limits_{\R^d\times \R^d}
\eins_{\sqrt{t}A \times \sqrt{t}A  }(\phi_t(x),\phi_t(y))m\otimes m(\dd x\dd y)\\
&=\limsup_{t \to \infty} \int\limits_{(\R^d\times \R^d) \setminus D} \P \big(\frac{1}{\sqrt t} (\phi_t(x),\phi_t(y)) \in A \times A\big) m\otimes m(\dd x\dd y) \le \mathcal{N}(0,\text{Id}_{\R^{2d}}) (A \times A) \\
&= \left( \mathcal{N}(0, \text{Id}_{\R^d}) (A)\right)^2
\end{align*}
 and hence
\begin{align*}
0 \le &\E(\left[l_t(A)- \mathcal{N}(0,\text{Id}_{\R^d}) (A)\right]^2)\\
&= \E\, l_t^2(A)- 2\mathcal{N}(0,\text{Id}_{\R^d}) (A)\E \, l_t(A) + \left( \mathcal{N}(0,\text{Id}_{\R^d}) (A)\right)^2 \underset{t \to +\infty}\longrightarrow 0\,.
\end{align*}
Taking complements, the assertion of the theorem follows also for open subsets. The case of a general Borel subset
$A$ follows since every probability measure on $\R^d$ is regular. The ``in particular'' statement in
the Theorem now follows easily and the proof is complete.\hfill $ \square $\\[0.5cm]

\noindent\textbf{Proof of Theorem \ref{an_main2}}:\\
IBFs enjoy the property that their laws are invariant under time reversal. This fact can be easily derived (see e.g. \cite{Dim06}, Corollary 1.2.1) by observing that the correction term in the backwards generator vanishes (\cite{Ku90}, Section 4.2). In particular, for each fixed $t>0$, the random functions $\phi_t$ and $\phi_t^{-1}$ have the same distribution and therefore
\begin{align} \label{an_main2_1}
&\E \left[\lambda(\phi_t(B) \cap \sqrt{t}A )\right ]^2=\E \int \eins_{\phi_t(B)\times \phi_t(B)}(x,y) \eins_{\sqrt{t}A\times \sqrt{t}A}(x,y) \dd x\dd y \nonumber \\[0.4cm]
&= \E \int \eins_{B\times B}(\phi_t(x),\phi_t(y)) \eins_{\sqrt{t}A\times \sqrt{t}A}(x,y) \dd x\dd y\,.
\end{align}
Let $\mu(\dd x \dd y):=\psi(|x-y|)\dd x\dd y$ be the  invariant measure for the two-point motion given in (\ref{ibf2p_inv_meas}). The reversibility of the two-point motion with respect to $\mu(\dd x \dd y)$  implies
\begin{align} \label{an_main2_2}
& \E \int \eins_{B\times B}(\phi_t(x),\phi_t(y)) \eins_{\sqrt{t}A\times \sqrt{t}A}(x,y)  \dd x\dd y \nonumber \\[0.4cm]
&= \E \int \eins_{B\times B}(\phi_t(x),\phi_t(y)) \eins_{\sqrt{t}A\times \sqrt{t}A}(x,y)\frac{1}{\psi(|x-y|)} \mu(\dd x\dd y) \nonumber \\[0.4cm]
&= \int \eins_{B\times B} (x,y)
\E \left\{\eins_{\sqrt{t}A\times \sqrt{t}A}(\phi_t(x),\phi_t(y))\frac{1}{\psi(|\phi_t(x)-\phi_t(y)|)}\right\} \mu(\dd x\dd y)\,.
\end{align}
Setting
$$M_t(A):=\eins_{\sqrt{t}A\times \sqrt{t}A}(\phi_t(x),\phi_t(y))\frac{1}{\psi(|\phi_t(x)-\phi_t(y)|)}$$ and
$$D_L:=\{(x,y)\in \R^{2d} : |x-y|\ge L,\, |x|\ge L,\, |y|\ge L\},\;\; L>0$$ we have
\begin{align} \label{an_main2_3}
& \E M_t(A) = \E (M_t(A)\eins_{D_L^c}(  \phi_t(x),\phi_t(y) ))+\E (M_t(A)\eins_{D_L}(  \phi_t(x),\phi_t(y) )).
\end{align}
 The first term in the decomposition converges to zero by the Lebesgue's bounded convergence theorem,
since $\frac{1}{\psi(|x-y|)}$ is uniformly bounded from above (see  \eqref{betta} and \eqref{ibf2p_inv_meas_asymp}) and
$\eins_{D_L^c}(  \phi_t(x),\phi_t(y) )\to 0$ in probability as $t \to +\infty$. That is,
\begin{align} \label{an_main2_4}
\lim_{t \to +\infty}\E (M_t(A)\eins_{D_L^c}(  \phi_t(x),\phi_t(y) ))=0\,.
\end{align}
The second term in (\ref{an_main2_3}) can be bounded from  above by:
\begin{align} \label{an_main2_5}
\begin{split}
& H_L \E \eins_{\sqrt{t}A\times \sqrt{t}A}(\phi_t(x) ,\phi_t(y)),
\end{split}
\end{align}
where $H_L:=\sup_{(x,y)\in D_L}\big \{\frac{1}{\psi(|x-y|)}\big \}$.
Now we assume that the set $A$ is closed. According to Lemma \ref{an_CLTforIBF}, we have
\begin{align} \label{an_main2_6}
& \lim_{t \to \infty} H_L \E \eins_{\sqrt{t}A\times \sqrt{t}A}(\phi_t(x) ,\phi_t(y)) \le H_L\left( \mathcal{N}(0,\text{Id}_{\R^{d}})(A)\right)^2\,.
\end{align}

Combining (\ref{an_main2_6}), (\ref{an_main2_4}) and (\ref{an_main2_3}), we obtain for arbitrary $L > 0$
\begin{align} \label{an_main2_8}
\limsup_{t \to +\infty}  \E M_t(A) \le H_L\left( \mathcal{N}(0, \text{Id}_{\R^{d}})(A)\right)^2.
\end{align}
We have $H_L\searrow 1$ for $L\to +\infty$ because $\psi(s)\to 1$ as $s \to \infty$. Since (\ref{an_main2_8}) holds for every $L>0$  we can conclude that
\begin{align} \label{an_main2_10}
&\limsup_{t \to +\infty}  \E \eins_{\sqrt{t}A\times \sqrt{t}A}(\phi_t(x) ,\phi_t(y))
\frac{1}{\psi(|\phi_t(x)-\phi_t(y)|)} \le \left( \mathcal{N}(0,\text{Id}_{\R^{d}})(A)\right)^2\,.
\end{align}
Plugging  this into (\ref{an_main2_2}) and applying  Fatou's lemma yields
\begin{align} \label{an_main2_11}
\begin{split}
 \limsup_{t \to +\infty} \E \left(\lambda(\phi_t(B) \cap \sqrt{t}A )\right )^2&= \limsup_{t \to +\infty}  \int \eins_{B\times B} (x,y)\E M_t(A) \mu (\dd x\dd y) \\
&\le  \mu(B\times B)\left( \mathcal{N}(0,\text{Id}_{\R^{d}})(A)\right)^2.
\end{split}
\end{align}
The application  of Fatou's lemma is justified since  $1/\psi$ is bounded (and therefore $M_t(A)$ is uniformly bounded) and $\mu(B \times B)<\infty$. Indeed,
$$
\mu(B \times B)=\int\limits_{\R^d\times \R^d}\eins_B(x)\eins_B(y)\psi(|x-y|)\dd x\dd y=\int\limits_{\R^d\times \R^d}\eins_B(x)\eins_B(z+x)\psi(|z|)\dd x\dd z <\infty,
$$
since $\lambda(B)<\infty$ and $\psi(|z|)$ is integrable over the unit $d$-dimensional ball B$(0,1)$ and bounded outside. Notice, that the condition $\lambda_1>0$ is really needed at this point:  clearly

$$
\int\limits_{\text{B}(0,1)}\psi(|z|)\dd z =\frac{2\pi^\frac{d}{2}}{\Gamma(\frac{d}{2})}\int\limits_0^1 \psi(x)x^{d-1}\dd x
$$
and according to  (\ref{ibf2p_inv_meas_asymp}) for $x \to 0$  we have $\psi(x)x^{d-1} \sim \frac{c}{\beta_L}x^{(d-1)\frac{\beta_N}{\beta_L}-2}$. The condition that $\lambda_1=(d-1)\beta_N/2-\beta_L/2>0$ clearly reads  $(d-1)\frac{\beta_N}{\beta_L}-2>-1$ and thus implies that $\psi(x)x^{d-1}$ is integrable in a neighborhood of the origin. The above reasoning also  implies that the condition $\lambda_1>0$ is necessary for  $\psi$ to be  integrable over the unit $d$-dimensional ball B$(0,1)$.
\\

As in the derivation of (\ref{an_main2_11}), we  get for an arbitrary (not necessarily closed) Borel set $A$:
\begin{align} \label{an_main2_13}
&\E \big[\mathcal{N}(0,\text{Id}_{\R^d})(A) \lambda( \phi_t(B))\lambda(\phi_t(B)\cap \sqrt{t}A)\big] \nonumber \\[0.4cm]
&=\mathcal{N}(0,\text{Id}_{\R^d})(A)  \E \int \eins_{\phi_t(B)\times \phi_t(B)}(x,y) \eins_{\sqrt{t}A}(x) \frac{1}{\psi(|x-y|)}\mu(\dd x\dd y) \nonumber \\[0.4cm]
&=\mathcal{N}(0, \text{Id}_{\R^d})(A) \int \eins_{B\times B}(x,y) \E \left[ \eins_{\sqrt{t}A}(\phi_t(x)) \frac{1}{\psi(|\phi_t(x)-\phi_t(y)|)}\right]\mu(\dd x\dd y) \nonumber \\[0.4cm]
 & \underset{t \to +\infty} \longrightarrow (\mathcal{N}(0,\text{Id}_{\R^d})(A))^2 \mu(B\times B) \,.
\end{align}
Combining (\ref{an_main2_11}), (\ref{an_main2_13}) and Theorem \ref{ibf_vol_thm0}(ii), we obtain for a closed set $A$
\begin{align*}
0 \le &\limsup_{t \to \infty} \E\left[  \lambda( \phi_t(B)\cap \sqrt{t}A )-\mathcal{N}(0,\text{Id}_{\R^d})(A)\lambda(\phi_t(B)) \right]^2
=\limsup_{t \to \infty} \Big\{ \E\left[  \lambda( \phi_t(B)\cap \sqrt{t}A )\right]^2\\[0.4cm]
&+ \E\left[ \mathcal{N}(0,\text{Id}_{\R^d})(A)\lambda(\phi_t(B)) \right]^2-2\E\left[  \lambda( \phi_t(B)\cap \sqrt{t}A )\mathcal{N}(0,\text{Id}_{\R^d})(A)\lambda(\phi_t(B)) \right]\Big\}\\[0.4cm]
& \le 2(\mathcal{N}(0, \text{Id}_{\R^d})(A))^2 \mu(B\times B)-2(\mathcal{N}(0, \text{Id}_{\R^d})(A))^2 \mu(B\times B)=0\,,
\end{align*}
which proves (\ref{qares1}) for closed $A$. Taking complements, (\ref{qares1}) follows also for open sets and --
like in the proof of the previous theorem --
$\mathrm{L}^2$-convergence for a general Borel set follows from the regularity of probability measures on a Euclidean space.
Further (\ref{qares2}) is a consequence of (\ref{qares1}) and the fact that $\lambda(\phi_t(B))$ converges almost surely to a
nonnegative random variable (by Theorem \ref{ibf_vol_thm0}(ii)).\hfill $ \square $\\

\section{Persistence of volume  }
We investigate the evolution of  volume under the action of an IBF, that is we are interested in the asymptotics of
$$V_t(A):=\int \limits_{\R^d}\eins_{\{\phi_t(A)\}}(x)\,\, \lambda (\dd x)=\lambda(\phi_t(A))\,,$$
for $A\in \B(\R^d)$. Note that $V_t$ defines a random measure on the Borel subsets of $\R^d$.

We aim at showing $\P(\lim_{t\to \infty} V_t(A)>0)=1$ for a nonempty open set $A$ under additional conditions on the parameters of the IBF.
\subsection{The statistical equilibrium}

To put our result in the framework of existing results, in the following theorem we  provide  a collection of facts concerning  $\left(V_t\right)_{t \ge 0}$. They can be found in  works of Baxendale and Harris \cite{BaH86}, Baxendale \cite{Bax91},
Le Jan \cite{LeJ85}, \cite{LeJ87}, Darling and Le Jan \cite{DaLe88} and  Kunita \cite{Ku90} (Lemma 4.3.1 and Theorem 4.3.6),
where also the proofs are available.  The proofs can be found also in \cite{Dim06}, Theorem 4.0.1 -- with the
exception of the of part (iv) (which is contained in  \cite{DaLe88}).

\begin{theo}\label{ibf_vol_thm0}$ $
\begin{enumerate}
\item[\textnormal{(i)}] The family $\{V_t : t \ge 0\}$ converges a.s.  vaguely to a random measure $V$ as  $ t\to +\infty $,
which we call the {\em forward statistical equilibrium}.
\item[\textnormal{(ii)}]  Assume further that the top Lyapunov exponent $\lambda_1$  of the flow is strictly positive. Then for all bounded and measurable $f,g :\R^d \to  \R$
\begin{align*}
&\E\int\limits_{\R^d} f(x) \,\,V(\dd x) = \int\limits_{\R^d} f(x) \,\,\lambda(\dd x) \text{ and } \\
&\E \int\limits_{\R^{2d}} f(x)g(y) \,\,V \otimes V(\dd x\dd y) = \int\limits_{\R^{2d}} f(x)g(y)\psi(|x-y|) \,\,\lambda(\dd x \dd y)\,\,,
\end{align*}
where $\psi(|x-y|) \,\,\lambda(\dd x)$ is the invariant measure for the two-point motion \\$\big\{(\phi_t(x),\phi_t(y))^T : t \ge 0\}$, introduced in (\ref{ibf2p_inv_meas}). In other words, the first moment of
$V$ is the invariant measure for the one-point motion and the second moment of $V$ is the invariant measure for the two-point motion.\\
Furthermore, for arbitrary $A\in \B(\R^d)$, $\lambda(A) < \infty$, the process $t \mapsto V_t(A)$ is an $\mathrm{L}^2$-bounded martingale and
$$ \lim_{t\to +\infty} V_t(A)=V(A) \,\,\,\,\,\P-\text{almost surely and in }\mathrm{L}^2. $$

\item[\textnormal{(iii)}] The measure $V$  is almost surely orthogonal to the Lebesgue measure, provided it is neither zero nor Lebesgue measure itself (volume preserving case).
\item[\textnormal{(iv)}] In case the top Lyapunov exponent is strictly negative, the measure $V$ is almost surely identically zero.
\end{enumerate}
\end{theo}

\begin{remark}
The statistical equilibrium of a Brownian flow with stationary increments is usually defined
(see \cite{Bax91}, \cite{LeJ87}, \cite{DaLe88}) as the almost sure vague limit
$$\mu_s:=\lim_{t\to +\infty} \mu\circ \phi_{-t,0}^{-1},$$
where $\mu$ is an invariant measure for the one point motion (modulo a multiplicative constant the Lebesgue measure in case $\phi$ is an IBF).
In order to define the above limit one considers the canonical extension of the flow to double sided time, using the
independence and the stationarity of the increments. Time reversibility of an IBF implies that the processes
$\lambda \circ \phi_{-t,0}^{-1}$  and $V_t$ have the same joint distribution and therefore the statistical equilibrium
$\mu_s$ has the same law as the forward statistical equilibrium $V$. Observe that $V$ and $\mu_s$ are independent, since they are
functions of the increments of $\phi$ on $[0,\infty)$ respectively $(-\infty,0]$.
\end{remark}

\subsection{Persistence of volume}\label{persistence}

Let $A$ be a nonempty bounded open subset of $\R^d$. We will investigate, if the volume of the set $A$ persists almost surely under
the action of the flow, that is if
\begin{align*}
\P\left(\lim_{t \to \infty}\lambda(\phi_t(A))\ne 0\right)=1 \,.
\end{align*}
The statement is trivially true in the volume preserving case, characterized via $\frac{\beta_N}{\beta_L}= \frac{d+1}{d-1}$.
Darling and Le Jan showed in  \cite{DaLe88}, that in case $\lambda_1<0$, we have $$ \P\left(\lim\limits_{t \to \infty}\lambda(\phi_t(A))= 0\right)=1. $$
For the rest of the section, we will assume that $\lambda_1 > 0$.
The idea of the proof of the persistence of volume is to use a second moment method together with some geometric considerations
about the image of the set $A$ under the action of the flow. In fact, for an arbitrary set $A$ of positive Lebesgue measure,
we have by Theorem \ref{ibf_vol_thm0}
\begin{align*}
\lambda(A)=\E\left[ V(A)\right ]=\E\left[ V(A)\eins_{\{V(A)\ne 0\}}\right ]\le \sqrt{\E\left[ V^2(A)\right ]\P(V(A)\ne 0)}\,
\end{align*}
and therefore
\begin{equation}\label{absch}
\P(V(A)\ne 0) \ge \lambda^2(A)\frac{1}{\E\left[ V^2(A)\right ]}=\lambda^2(A)\frac{1}{\int\limits_{A \times A}\psi(|x-y|)\dd x\dd y}\,.
\end{equation}
Conditioning on $\F_t$, one obtains that $\P(V(A)\ne 0)=1$ provided that
\begin{align*}
\frac{1}{\lambda^2(\phi_t(A))}\int\limits_{\phi_t(A)\times \phi_t(A) }\psi(|x-y|)\dd x\dd y \underset{t\to \infty}\longrightarrow 1 \,.
\end{align*}
Unfortunately, our knowledge about the shape of $\phi_t(A)$ is not detailed enough to prove this in general. Let $\text{B}(z,r)$
denote the closed ball of radius $r$ and center $z$. Unless the flow preserves Lebesgue measure,
we have
\begin{align*}
\frac{1}{\lambda^2(\text{B}(z,r))}\int\limits_{\text{B}(z,r)\times \text{B}(z,r)}\psi(|x-y|)\dd x\dd y
\overset{r \downarrow 0}\longrightarrow +\infty\,,
\end{align*}
since $\psi(s)\to +\infty$ as $s \to +0$, so it will not help to replace $\phi_t(A)$ by a (possibly very small) ball contained in $\phi_t(A)$
in case $A$ is nonempty and open. Replacing $\phi_t(A)$ by long and thin tubes contained in  $\phi_t(A)$ will work however in some cases.\\
To prove the main theorem in this section, we first need the following two lemmas:\\
\begin{lemma}\label{ibf_vol_lemma1}
Let $h\colon \R_+ \to \R_+$ be a  monotonically decreasing function such that there are positive constants $a,C>0$ and a number $\mu>-1$ such that
$$ \lim_{s \to +\infty}h(s)=1 \,\,\text{ and }\,\,h(s)\le C s^\mu \,\text{ for }\, s\in(0,a) \,.$$
Let further   $Z:=Z(L,\delta)$ be the cylinder
$$Z(L,\delta):=\{x\in \R^d \, \colon \, x_1 \in [-L/2,L/2], \,|(0,x_2,\dots,x_d)^T|\le \delta\},\;\;L,\delta >0.$$
Then
\begin{enumerate}
\item[\textnormal{(i)}] $\lim\limits_{L\to +\infty}\frac{1}{2L}\int \limits_{-L}^L h(|r|)\dd r=1 $ and
\item[\textnormal{(ii)}] $\frac{1}{\lambda^2(Z(2L,\delta))}\int\limits_{Z(2L,\delta)\times Z(2L,\delta)}h(|x-y|)\dd x\dd y\le \frac{1}{2L}\int \limits_{-L}^L h(|r|)\dd r$.
\end{enumerate}
\end{lemma}

\begin{remark} \label{ibf_vol_rem1}
Clearly \textnormal{(i)} and \textnormal{(ii)} imply
\begin{align}\label{pers_vol_1}
\frac{1}{\lambda^2(Z(2L,\delta))}\int\limits_{Z(2L,\delta)\times Z(2L,\delta)}h(|x-y|)\dd x\dd y \underset{L \uparrow +\infty}\longrightarrow 1 \, \,\,\text{ uniformly in  } \delta\,,
\end{align}
since $h \ge 1$.

The claim in (\ref{pers_vol_1}) is   invariant with respect to rigid motions on $ \R^d $ and therefore holds for arbitrary cylinders with length $2L$ and radius $\delta$.
\end{remark}
\noindent \textbf{Proof of Lemma \ref{ibf_vol_lemma1}}:\\
(i) is obvious.
To show (ii), observe that
\begin{align*}
&\frac{1}{\lambda^2(Z(2L,\delta))}\int\limits_{Z(2L,\delta)\times Z(2L,\delta)}h(|x-y|)\dd x\dd y=\frac{1}{\lambda^2(Z(2L,\delta))}\int\limits_{Z(2L,\delta)}\int \limits_{Z(2L,\delta)-y }h(|x|)\dd x\dd y \nonumber \\[0.5cm]
&\le \frac{1}{\lambda(Z(2L,\delta))}\sup_{y \in Z(2L,\delta)}\int \limits_{Z(2L,\delta)-y }h(|x|)\dd x = \frac{1}{\lambda(Z(2L,\delta))}\int \limits_{Z(2L,\delta)}h(|x|)\dd x \nonumber \\
&\le \frac{1}{\lambda(Z(2L,\delta))}\int \limits_{Z(2L,\delta)}h(|x_1|)\dd x = \frac{1}{2L} \int_{-L}^{L} h(|r|) \dd r\, ,
\end{align*}
where we used the monotonicity of $h$ and the fact that $Z(2L,\delta)$ is centered at the origin.
\hfill $\square$\\[0.5cm]

\begin{lemma}\label{ibf_vol_lemma2}

Let the function $h$ be as in the previous lemma.  Let $A$ be a bounded subset of $\R^d$, diffeomorphic to an open ball and let $L$ be some
number strictly between 0 and the diameter of $A$.
Then there exist a positive integer $n$, $l>0$ and $n$ pairwise disjoint linear segments $\left(\gamma^i\right)_{i=1,\dots,n}$ of length $l$ and a number $\bar \delta>0$ such that:
\begin{itemize}
\item[(i)] the length of $\gamma(A):=\gamma:=\cup_{i=1}^n\gamma^i$ is at least $\frac{L}{7}$, that is $nl\ge \frac{L}{7}$,
\item[(ii)] for all $\delta \le \bar \delta$  the ``piecewise cylinder'' set $\gamma_\delta(A):=\gamma_\delta:=\cup_{i=1}^n\gamma^i_\delta$, where $\gamma^i_\delta$ is the closed cylinder with axis $\gamma^i$ and radius $\delta$ is completely contained in $A$. Moreover, the cylinders $ \gamma_\delta^i $ are disjoint.
\item[(iii)] the following holds for all $\delta \le \bar \delta$:
$$\frac{1}{\lambda^2(\gamma_\delta)}\int\limits_{\gamma_\delta\times \gamma_\delta}h(|x-y|)\dd x\dd y  \le  \frac{1}{\lambda^2(Z(nl,\delta))}\int\limits_{Z(nl,\delta)\times Z(nl,\delta)}h(|x-y|)\dd x\dd y\,.$$
\end{itemize}
\end{lemma}

\noindent\textbf{Proof}:\\
Since $A$ is an open connected subset of $\R^d$ with diameter greater than $L$, we can find a piecewise linear, connected curve
\begin{align*}
\tilde \gamma:=\bigcup_{i=1}^{\tilde n} \tilde{\gamma}^i
\end{align*}
with the properties:
\begin{enumerate}
\item[(1.)] The curve $\tilde \gamma$ is completely contained in $A$.
\item[(2.)] All linear segments $ \tilde{\gamma}^i $ have the same length $ l>0$, which is a strictly positive number depending on the shape of $A$.
\item[(3.)] The diameter of the curve  $\tilde \gamma$ is $L$.
\end{enumerate}

Let $ a, b\in \tilde \gamma$ are such that $|a-b|=L$ and  $Y_{ab}$ be the linear segment connecting $ a$ and $ b$, i.e.
$$Y_{ab}=\{x \in \R^d \,\colon \,x=\alpha  a + (1-\alpha)  b\,, \, \,\alpha \in [0,1]\}\,.$$
In order to ease notation and without loss of generality we can assume that $L=2n\cdot 3l$ for some integer $n$.
Consider a partition $\{ x_k : k=0,\dots,2n \}$ of $Y_{ab}$ by the points
$$x_k= a+k\frac{3l}{L}( b- a)=a+\frac{k}{2n}(b-a)\,,$$
that is $|x_{k+1}-x_k|=\frac{|b-a|}{2n}=3l$.\\
Let $E_k$ denote the $d-1$-dimensional  hyperplane containing $x_k$ and being orthogonal to $Y_{ab}$ and for all $k=1,\dots,2n$ let $V_k $ be the set of all points  strictly between the planes $E_{k-1}$ and $E_k$. For every $k \in \{1, \dots, 2n\}$ there is at least one  $i_k \in \{1, \dots ,\tilde n\}$ such that the line segment $\tilde{\gamma}^{i_k}$ is contained in $V_k$. Further, all these indices are different since the $V_k$'s are disjoint.\\
  Consider the disjoint union of  linear segments $\gamma \subset \tilde \gamma$ given by
\begin{align*}
\gamma=\bigcup_{k=1}^{n}\gamma^k \,\,\text{ with }\,\, \gamma^k=\tilde{\gamma}^{i_{2k-1}}\,,
\end{align*}
that is  $\gamma$ contains exactly one linear piece of length $l$ in every set $V_k$ for odd indices $k$, i.e. one linear segment in every second  $V_k$.\\
We  show that $\gamma $ satisfies (i), (ii) and (iii).\\
The different pieces in $\gamma$ are clearly disjoint. The length of $\gamma$ is given by
$$\text{Length}(\gamma)=nl=\frac{L}{6l}l=\frac{L}{6}>\frac{L}{7}\,.$$
Now,  the family $\left (\gamma^k \right)_{k=1,\dots,n}$ is disjoint and  contained in $A$ and therefore we can find $\delta_1>0$ such that for all  $\delta< \delta_1$
\begin{enumerate}
\item[1.] the family $\left (\gamma^k_\delta \right)_{k=1,\dots,n}$ is disjoint,
\item[2.] $\gamma_\delta:=\cup_{k=1}^n \gamma^k_\delta$ is contained in $A$,
\item[3.] for all $k \in \{1,\dots,n\}$ $\gamma^k_\delta$ is contained in  $V_{2k-1}$.
\end{enumerate}
It remains to show that (iii) holds. Intuitively this is clear, since if we piece together all cylinders $\gamma_\delta^k$ in order to
obtain a tube of length $nl$ and radius $\delta$, then we either reduce the distances between the points $x$ and $y$ if they  lie in different pieces
or  $|x-y|$ does not change if  $x$ and $y$ lie in the same piece.  Since $h$ is decreasing, (iii) follows. Those readers,
who really want to see a detailed rigorous proof of this, are referred to the proof of Lemma 4.1.2 in \cite{Dim06}.
\hfill $\square$\\[0.5cm]
Now we have collected all prerequisites for proving the main theorem in this section:

\begin{theo}\label{ibf_vol_thm1}
Let $\phi$ be an isotropic Brownian flow with
\begin{align*}
\frac{\beta_N}{\beta_L}>\frac{d}{d-1}\,.
\end{align*}
The volume of every nonempty open set $A$ persists with probability one under the action of the flow, i.e.
\begin{align*}
\P\big ( \,\lim_{t \to \infty}\lambda \left(\phi_t(A)\right)\ne 0\big)=1.
\end{align*}
\end{theo}
\begin{remark}\label{ibf_vol_rem3}
The parameters $\beta_L$ and $\beta_N$ always fulfill the inequalities
\begin{align*}
\frac{1}{3}\le \frac{\beta_N}{\beta_L}\le \frac{d+1}{d-1}\,,
\end{align*}
where $\beta_N(d-1)=\beta_L(d+1)$ corresponds to the volume preserving case, for which the persistence of the volume is trivial. The preceding theorem states that the volume persists also if we are close to the volume preserving case, in the sense that
\begin{align*}
\frac{d}{d-1}< \frac{\beta_N}{\beta_L}\le \frac{d+1}{d-1}\,.
\end{align*}
Recall that the top Lyapunov exponent is strictly positive iff
\begin{align*}
\frac{1}{d-1}< \frac{\beta_N}{\beta_L}\,,
\end{align*}
which is satisfied under the assumption of Theorem \ref{ibf_vol_thm1}.
\end{remark}

\noindent \textbf{Proof of Theorem \ref{ibf_vol_thm1}}:\\
Since every nonempty open set $A$ contains an open ball, it is enough to prove that
\begin{align*}
\P\big(\,\lim_{t \to \infty}\lambda \left(\phi_t(\text{B}(0,r))\right)\ne 0\big)=1 \text{ for every } r>0\,.
\end{align*}
Since $\frac{\beta_N}{\beta_L}>\frac{d}{d-1}$, we have $\psi(s) \sim \frac{c}{\beta_L} s^{\mu}$ for $s \downarrow 0$
with $\mu > -1$. Therefore there exists a function $h \ge \psi$ which satisfies all assumptions in Lemma
\ref{ibf_vol_lemma1} (with the same $\mu$). We abbreviate $\Gamma_s:=\gamma_{\delta}(\phi_s(\text{B}(0,r)))$,
and $Z_s:= Z(nl,\delta)$, where $\gamma_{\delta}$ and $Z(nl,\delta)$ are defined as in Lemma \ref{ibf_vol_lemma2}
in case we replace $A$ by $\phi_s(\text{B}(0,r))$.

Using the independence of the increments  of the flow, Lemmas \ref{ibf_vol_lemma1} and \ref{ibf_vol_lemma2},
and (\ref{absch}), we obtain:
\begin{align}\label{pers_vol_13}
&\P\big ( \,\lim_{t \to \infty} \lambda \left(\phi_{0,t}(\text{B}(0,r))\right) \ne 0 \,\,\big | \,\,\F_s\big)=\P\big (\, \lim_{t \to \infty} \lambda \left( \phi_{s,t}\circ \phi_{0,s}(\text{B}(0,r))\right)\ne 0 \,\,\big | \,\,\F_s\big)\nonumber \\[0.4cm]
&=\P\big ( \,\lim_{t \to \infty} \lambda \left(\phi_{s,t}(B)\right)\ne 0 \big )\Big |_{B=\phi_{0,s}(\text{B}(0,r))}\ge \P\big (\, \lim_{t \to \infty} \lambda \left( \phi_{s,t}(B)\right)\ne 0 \big )\Big |_{B=Z_s}\nonumber\\
&\ge \lambda^2 \left(Z_s \right) \left[ \int_{Z_s \times Z_s} \psi(|x-y|) \dd x \dd y\right]^{-1}
\ge \lambda^2 \left( Z_s \right) \left[\int_{Z_s \times Z_s} h(|x-y|) \dd x \dd y \right]^{-1}\underset{s \to \infty}
\longrightarrow 1,\nonumber
\end{align}
where we used the fact that the diameter of $\phi_s(\text{B}(0,r))$ converges to $\infty$ as $s \to \infty$ almost
surely (see \cite{CSS99}). Unconditioning, the assertion follows. \hfill $\square$

\bibliographystyle{abbrv}

\bibliography{biblio}

\end{document}